\documentclass[reqno,11pt]{amsart}
\usepackage{amsmath, amssymb, amsthm, mathrsfs}
\usepackage{graphicx,url}
\usepackage[all]{xy}

\newtheorem{defn}{Definition}
\newtheorem{theorem}{Theorem}

\newtheorem{lemma}{Lemma}
\newtheorem{example}{Example}

\newcommand{\ZZ}{\mathbb{Z}}
\newcommand{\QQ}{\mathbb{Q}}
\newcommand{\RR}{\mathbb{R}}
\newcommand{\CC}{\mathbb{C}}
\newcommand{\HH}{\mathbb{H}}
\newcommand{\SLZ}{\operatorname{SL}_2(\ZZ)}

\newcommand{\PSLZ}{\operatorname{PSL}_2(\ZZ)}

\newcommand{\im}{\operatorname {Im}}

\begin{document}

\title[Computations with finite index subgroups of $\text{PSL}_2(\ZZ)$]{Computations with finite index subgroups of $\text{PSL}_2(\ZZ)$ using Farey Symbols}
\author{Chris A. Kurth}
\address{Department of Mathematics\\Iowa State University\\Ames, IA 50011 \\USA}
\email{kurthc@iastate.edu}
\author{Ling Long}

\address{Department of Mathematics\\Iowa State University\\Ames, IA 50011 \\USA}
\email{linglong@iastate.edu}

\begin{abstract}Finite index subgroups of the modular group are of
great arithmetic importance.  Farey symbols, introduced by Ravi
Kulkarni in 1991, are a tool for working with these groups. Given
such a group $\Gamma$, a Farey symbol for $\Gamma$ is a certain
finite sequence of rational numbers (representing vertices of a
fundamental domain of $\Gamma$) together with pairing information
for the edges between the vertices. They are a compact way of
encoding the information about the group and they provide a simple
way to do calculations with the group. For example: calculating an
independent set of generators and decomposing group elements into a
word in these generators, finding coset representatives, elliptic
points, and genus of the group, testing if the group is congruence,
etc. In this expository article, we will discuss Farey Symbols and
explicit algorithms for
working with them. %The first author has implemented a collection of
%these functions for finite index subgroups of the modular group into
%a SAGE package which is
\end{abstract}

\maketitle

\section{Introduction} Modular forms are certain functions defined
on the upper half plane displaying certain symmetries under the
M\"{o}bius transformation action of a finite-index subgroup of
$\PSLZ$. The theory of modular forms has been in the central stage
of number theory for more than one century and continues to be one
of its most exciting areas. Working with modular forms requires
knowing information about their underlying groups. There is a vast
literature about many aspects of finite index subgroups of the
modular groups and their relations with other fields such as
combinatorics, algebraic curves. Interested readers are referred to
articles like \cite{a-sd}, \cite{birchb94}, or a recent survey
article by the second author on these groups and their modular forms
\cite{long06-survey}.

Some finite index subgroups of the modular group can be described
purely by congruence relations, and as such are called congruence
subgroups of the modular group. These groups are relatively easy to
work with, as they contain a certain normal subgroup $\Gamma(N)$,
such that the quotient of the group by $\Gamma(N)$ is just a
subgroup of $\text{PSL}_2(\ZZ/N\ZZ)$. Most computational methods for
working with modular forms work only for congruence subgroups (cf.
\cite{stein07}). Noncongruence subgroups are finite-index subgroups
of $\PSLZ$ that cannot be described by congruence relations. Even
though they make up the majority of finite-index subgroups of
$\PSLZ$ there are few tools for working with them. One tool that
works equally well with both congruence and noncongruence groups is
the method of Farey symbols introduced by Ravi Kulkarni
\cite{kulkarni91}. In this expository paper we will recast how to
use Farey symbols and some related computational topics. We will
discuss some explicit algorithms for working with Farey symbols. The
first author has implemented a collection of such algorithms into a
free SAGE package called ``KFarey". It should be  made clear to the
readers that ``KFarey" is an undergoing project and we will continue
to improve the current functions, implement other existing
algorithms such as \cite{hsu97, Lang02normofmg},  and  investigate
new algorithms to make ``KFarey"  useful to a wide range of
audience. The second author would like to thank  ICAC 2007
organizers for inviting her to attend the conference and give a talk
on this topic.

\section{Subgroups of $\PSLZ$} Let $\SLZ$ be the group of $2 \times
2$ matrices with integer coefficients and determinant 1, and let
$\PSLZ = \SLZ / \{I,-I\}$. Let $\HH$ be the upper half plane $\HH :=
\{z\in \CC : \im(z)>0 \}$. Then $\PSLZ$ acts faithfully on $\HH$
under the action
\begin{equation*}
\gamma z = \frac{az+b}{cz+d}
\end{equation*}
where $\gamma = \begin{pmatrix}a&b\\c&d \end{pmatrix} \in \PSLZ$.
The objects of our study will be the finite index subgroups of
$\PSLZ$. For example, the standard congruence groups, $\Gamma(N)$,
$\Gamma_1(N)$, and $\Gamma_0(N)$.

If $\Gamma$ is a finite index subgroup of $\PSLZ$ then the action of
$\Gamma$ partitions $\QQ \cup \{\infty\}$ into equivalence classes,
where $q_1 \sim q_2$ if $q_1 = \gamma q_2$ for some $\gamma \in
\Gamma$. These equivalence classes $\{q\}$ are called the
\textbf{cusps} of $\Gamma$ and the width of the cusp $\{q\}$ is
$[\text{Stab}_{\PSLZ} (q) : \text{Stab}_\Gamma (q)]$. We say that
the level of $\Gamma$ is the least common multiple of the cusp
widths of $\Gamma$.

Recall that a congruence subgroup of $\PSLZ$ is a subgroup $\Gamma$
that contains a principal level $N$ congruence subgroup $\Gamma(N)$
for some $N$. If $N$ is the smallest such $N$ such that this is
true, then $\Gamma$ has level $N$.

\begin{defn}
Let $\Gamma$ be a finite index subgroup of $\PSLZ$. For our
purposes, a fundamental domain of $\Gamma$ is a hyperbolic polygon
$P$ on $\HH \cup \QQ \cup \{\infty\}$ such that:
\begin{enumerate}
\item  If $z$ is in the interior of $P$ and $\gamma \in \Gamma$,
then $\gamma z \in P$ implies $\gamma = I$.
\item For every $z \in \HH$ there is $\gamma \in \Gamma$ such that
$\gamma z \in P$.
\end{enumerate}
\end{defn}

\begin{lemma}
Let $\Gamma$ be a finite index subgroup of $\PSLZ$ and $P$ a
hyperbolic polygon. Suppose $P$ is such that:
\begin{enumerate}
\item If $z$ is in the interior of $P$ and $\gamma \in \Gamma$,
then $\gamma z \in P$ implies $\gamma = I$.
\item For each side $e$ of $P$, there is $\gamma \in \Gamma$ such
that $\gamma$ maps $e$ to another side of $P$ in an orientation
reversing manner.
\end{enumerate}
Then $P$ is a fundamental domain of $\Gamma$.
\end{lemma}
\begin{proof}
Let the images $\gamma P$ of $P$ under elements $\gamma$ of $\Gamma$
be called $P$-tiles. By Condition 1 they cannot overlap. Also, given
Condition 1, we only need to show that $\displaystyle \HH \subseteq
\bigcup_{\gamma \in \Gamma} \gamma P$. Suppose this is not true.
Then there is $\eta \in \Gamma$ such that $\eta P$ has an edge $e$
without a $P$-tile on the other side. Then $\eta^{-1} e$ is an side
of $P$ and the $\gamma$ of Condition 2 maps $\eta^{-1} e$ to another
side of $P$. Specifically, $\gamma^{-1}$ maps $P$ to a $P$-tile
adjacent to $P$ across the side $e$. Then $\eta \gamma^{-1} P$ is a
$P$-tile adjacent to $\eta P$ across the side $e$. Contradiction.
\end{proof}

\section{Farey Symbols}
\subsection{Special Polygons}
Farey Symbols were introduced by Ravi Kulkarni in 1991
\cite{kulkarni91} as a compact and efficient way to compute with
finite index subgroups of $\PSLZ$. The idea is to describe the group
by a fundamental domain with vertices at certain rational numbers
and certain hyperbolic arcs joining these rational numbers. Most of
the theory here is summarized from \cite{kulkarni91}.

If $x$ and $y$ are two points on $\HH \cup \QQ$ then there is a
unique circle passing through $x$ and $y$ with center on $\QQ$. We
say the \textbf{hyperbolic arc} joining $x$ and $y$ is the arc of
this circle contained in $\HH \cup \QQ$ joining $x$ and $y$. We also
say the hyperbolic arc joining $x\in \HH$ to $\infty$ is the
vertical line segment $\{x+ti:0 \leq t\in \RR \} \cup \{\infty\}$.
We write $H_{x,y}$ for the hyperbolic arc joining $x$ and $y$. A
\textbf{hyperbolic polygon} is a polygon composed of hyperbolic
arcs.

%%% Mention that the image of a hyperbolic arc under \gamma is a hyperbolic arc.

Through the course of this paper, when a vertex of a hyperbolic arc
is in $\QQ$ it will always be assumed to be in the form
$\frac{a}{b}$ with $a,b\in \mathbb Z$, $(a,b) = 1$ and $b > 0$. If
the vertex is $\infty$, we will write it either $\frac{-1}{0}$ or
$\frac{1}{0}$ (depending on if it is the leftmost or rightmost
element of a Farey sequence).

Let $\rho = \frac{1}{2} + \frac{\sqrt{3}}{2} i$ and let $T$ be the
hyperbolic triangle with vertices $\rho$, $\rho^2$ and $\infty$.
Then $T$ is a fundamental domain for $\PSLZ$ (\cite{kob1} Prop.
III.1). Let $E_e$ be the edge joining $i$ to $\infty$, $E_o$ be the
edge joining $\rho$ to $\infty$, and $E_f$ be the edge joining $i$
to $\rho$. Then we call an arc $A$ in the upper half plane an
\textbf{even edge} (resp. \textbf{odd edge}, resp. \textbf{f-edge})
if $A = \gamma E_e$ (resp. $A = \gamma E_o$, resp. $A = \gamma E_f$)
for some $\gamma \in \PSLZ$ (See Figure 1). $E_e$ and
$\begin{pmatrix}0&-1\\1&0\end{pmatrix} E_e$ together form a
hyperbolic arc from $0$ to $\infty$, and in general even edges come
in pairs joining rational numbers $\frac{a}{b}$, and $\frac{a'}{b'}$
with $|a'b-ab'| = 1$ because of the following lemma:

\begin{figure}[t]
 \includegraphics*[angle=270,scale=.4, viewport=50 200 400 600]{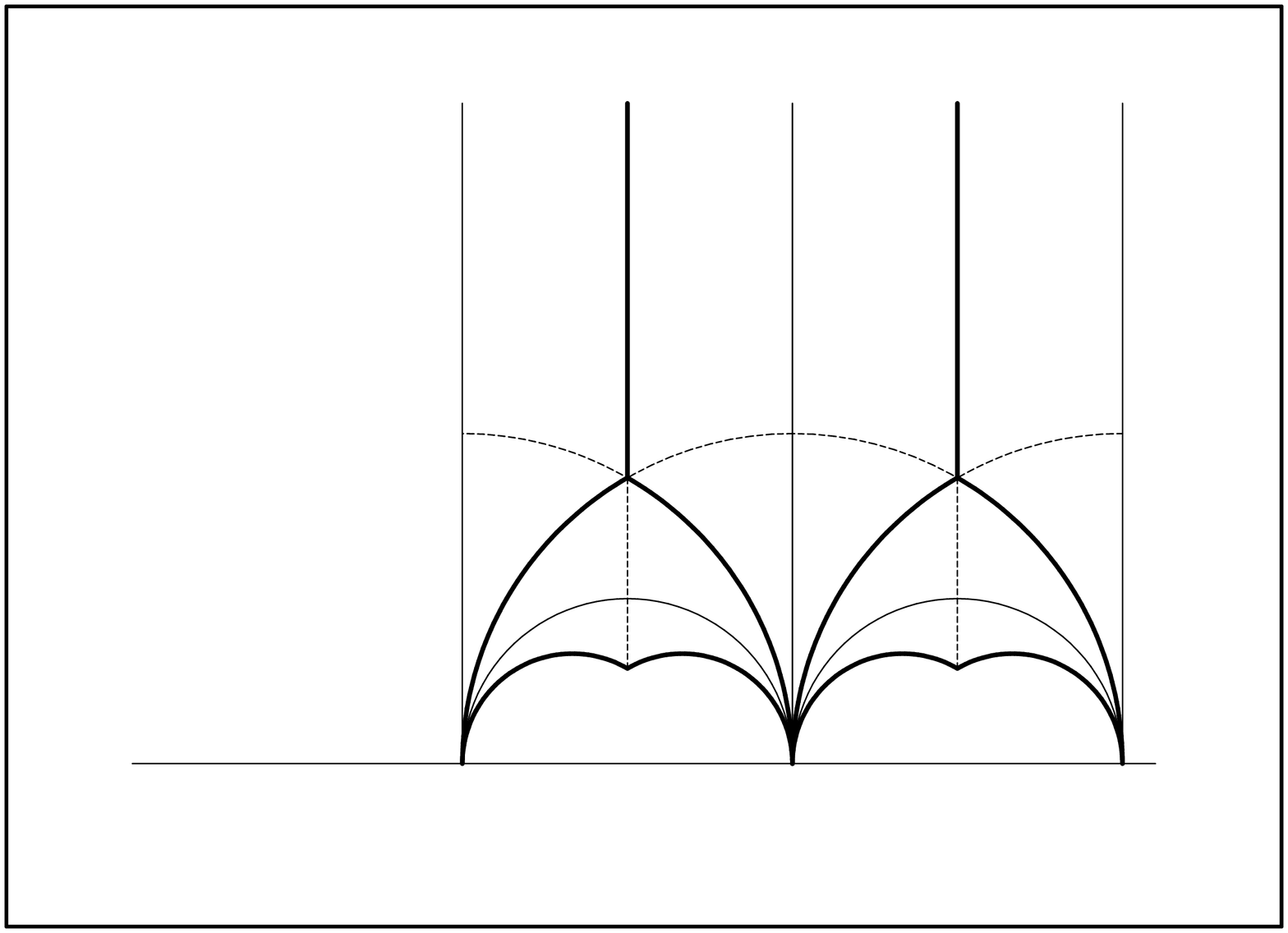}
 \caption{Even edges are thin, odd edges are thick, and f-edges are dashed}
\end{figure}

\begin{lemma}
If $\gamma \in \PSLZ$ and $a_1/b_1$, $a_2/b_2$, $a'_1/b'_1$, and
$a'_2/b'_2$ are rational numbers in simplest form such that
$$\gamma(a_1/b_1) = a'_1/b'_1, \text{ and } \gamma(a_2/b_2) = a'_2/b'_2$$ then
\begin{equation*}
a_2 b_1 - a_1 b_2 = a'_2 b'_1 - a'_1 b'_2
\end{equation*}
\end{lemma}
\begin{proof}
If $\gamma = \begin{pmatrix}A&B\\C&D\end{pmatrix}$ then
$$\gamma(\frac{a_1}{b_1}) = \frac{A a_1 + B b_1}{C a_1 + D b_1} =
\frac{a'_1}{b'_1} \text{ and } \gamma(\frac{a_2}{b_2}) = \frac{A a_2
+ B b_2}{C a_2 + D b_2} = \frac{a'_2}{b'_2}.$$ So:
\begin{eqnarray*}
a'_1 b'_2 - a'_2 b'_1 &=& (Aa_1 + Bb_1)(Ca_2 + Db_2) - (Aa_2 +
Bb_2)(Ca_1 + Db_1) \\
 &=& ADa_1b_2 + BCa_2b_1 - ADa_2b_1 - BCa_1b_2\\
 &=& (AD-BC)(a_1b_2 - a_2 b_1)\\
 &=& (a_1b_2 - a_2 b_1)
\end{eqnarray*}
\end{proof}

So the quantity $a_2 b_1 - a_1 b_2$ is invariant under
transformations in $\PSLZ$. Note that even edges, odd edges and free
edges only map to even edges, odd edges and free edges respectively
under transformations $\gamma \in \PSLZ$.

\begin{defn}
A special polygon $P$ is a convex hyperbolic polygon together with a
side pairing defined in the following way: The polygon is such that:
\begin{enumerate}
\item The boundary of $P$ consists of even and odd edges.
\item The even edges of $P$ come in pairs, each pair forming a
hyperbolic arc between elements of $\QQ \cup \{\infty\}$.
\item The odd edges of $P$ come in pairs, each pair meeting a vertex
with initial angle $\frac{2\pi}{3}$.
\end{enumerate}
The sides of the polygon are denoted as follows:
\begin{enumerate}
\item Each odd edge is called an \textbf{odd side}.
\item As even edges come in pairs, either each edge of the pair is
an \textbf{even side}, or the union of the two edges (a semicircle)
is called a \textbf{free side}.
\end{enumerate}
The side pairing on the edges is defined as follows:
\begin{enumerate}
\item Each odd side is paired with the odd side it meets at an
angle of $\frac{2\pi}3$. This is called an \textbf{odd pairing}.
\item Each even side is paired with the even side with which it
forms a semicircular arc. This is called an \textbf{even pairing}.
\item There are an even number of free sides and they are
partitioned into sets of two, each called a \textbf{free pairing}.
\end{enumerate}
\end{defn}

We will always assume that $0$ and $\infty$ are vertices of $P$.

The sides of a special polygon $P$ have a natural orientation
obtained by tracing the perimeter of the polygon in a certain
direction. If $\{s,s'\}$ is a side pairing then there is a unique
$\gamma \in \PSLZ$ such that $\gamma$ maps $s$ to $s'$ in an
orientation-reversing manner. We call this the \textbf{side pairing
transformation} associated with the side pairing, and we let
$\Gamma_P$ be the group generated by all the side pairing
transformations of $P$. Note that it doesn't matter which side we
pick for $s$ and which for $s'$ because the two possible $\gamma$'s
are inverses of each other. Also note that if $s$ is an even side
(resp. odd side) then $\gamma$ is order 2 (resp. order 3).

Two theorems of Kulkarni are fundamental here:

\begin{theorem}[\cite{kulkarni91} Theorem 3.2]
If $P$ is a special polygon then $P$ is a fundamental domain for
$\Gamma_P$. Moreover, the side pairing transformations $\{\gamma_i
\}$ are an independent set of generators of $\Gamma_P$ (i.e. the
only relations on the $\gamma_i$'s are $\gamma_i^2 = 1$ or
$\gamma_i^3 = 1$ for any finite-order $\gamma_i$'s).
\end{theorem}

\begin{theorem}[\cite{kulkarni91} Theorem 3.3]
For every $\Gamma \subset \PSLZ$ of finite index, there is a special
polygon $P$ such that $\Gamma = \Gamma_P$.
\end{theorem}
\begin{proof}\cite{kulkarni91} and also follows from the proof of
the algorithm in Section 4.
\end{proof}

Note that although it is true that any subgroup of $\PSLZ$ with
fundamental domain $F$ is generated by the transformations that map
its edges together, the fact that the set of generators of a special
polygon is an independent set of generators is something special to
the special polygon. For example $\Gamma(2)$ has a fundamental
domain shown in Figure 2.
\begin{figure}
 \includegraphics*[angle=270,scale=.4, viewport=100 130 450 530]{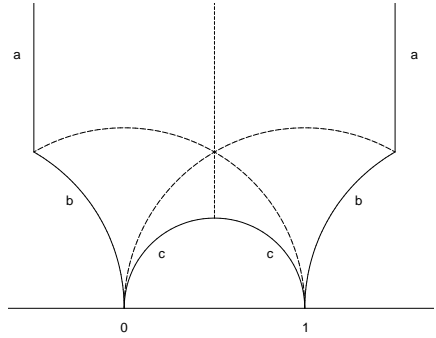}
 \caption{A fundamental
domain for $\Gamma(2)$}
\end{figure}
There are six sides, and the three side pairing transformations are
$\begin{pmatrix}1&2\\0&1\end{pmatrix}$,
$\begin{pmatrix}1&0\\2&1\end{pmatrix}$, and
$\begin{pmatrix}3&-2\\2&-1\end{pmatrix}$ But this is not a
independent list of generators because
$\begin{pmatrix}3&-2\\2&-1\end{pmatrix}^{-1}
\begin{pmatrix}1&2\\0&1\end{pmatrix} =
\begin{pmatrix}1&0\\2&1\end{pmatrix}$.

A special polygon for $\Gamma(2)$ is shown in Figure 3.
\begin{figure}
 \includegraphics*[angle=270,scale=.4, viewport=100 200 450 630]{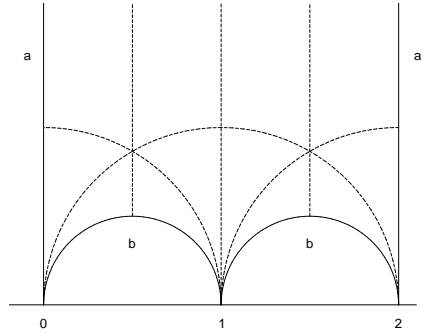}
 \caption{ A special polygon for
$\Gamma(2)$} \end{figure}
The pairing transformations from the
special polygon are $\begin{pmatrix}1&2\\0&1\end{pmatrix}$ and
$\begin{pmatrix}3&-2\\2&-1\end{pmatrix}$. These are independent
generators of $\Gamma(2)$.

\subsection{Farey Symbols}
Recall that the classical Farey sequences $F_n$ are constructed by
taking all the rational numbers $0 \leq a/b \leq 1$ with denominator
at most $n$ and $(a,b) = 1$ and writing them as a finite sequence in
ascending order $\{a_0/b_0, \dots, a_n/b_n \}$. Then for each $i$ we
have $a_{i+1}b_i - a_i b_{i+1} = 1$. We are interested in sequences
that satisfy this condition.

\begin{defn}
A \textbf{generalized Farey sequence} is a finite sequence:
$$\{\frac{-1}{0}, x_0, \dots, x_n, \frac{1}{0} \}$$ such that:
\begin{enumerate}
\item Each $x_i = a_i/b_i$ is a rational number in reduced form with
$b_i>0$. Additionally, we often consider $x_{-1} = \frac{-1}{0}$ and
$x_{n+1} = \frac{1}{0}$.
\item If we let  $a_{-1} = -1$, $b_{-1} = 0$, $a_{n+1} = 1$, and $b_{n+1} =
0$ then
\begin{equation}
a_{i+1}b_i - a_i b_{i+1} = 1
\end{equation}
for $-1 \leq i \leq n$.
\end{enumerate}
\end{defn}

Note that this definition forces $x_0$ and $x_n$ to be integers. We
will always assume $x_i = 0$ for some $i$.

\begin{defn}
A \textbf{Farey symbol} is a generalized Farey sequence with some
additional pairing information. Namely, between each adjacent
entries $x_{i-1}$ and $x_i$ we assign a \textbf{pairing} $p_i$ which
is either a positive integer called a \textbf{free pairing} or the
symbol ``$\circ$'' called an \textbf{even pairing} or ``$\bullet$''
called an \textbf{odd pairing}. Each integer that appears as a free
pairing appears exactly twice in the pairing information.
\end{defn}

So if $P$ is a special polygon, let $x_0, \dots, x_n$ be the
vertices of $P$ lying in $\QQ$ listed in ascending order. Recall
these vertices satisfy $a_{i+1}b_i - a_i b_{i+1} = 1$. Then
$\{\frac{-1}{0}, x_0, \dots, x_n, \frac{1}{0} \}$ is a generalized
Farey sequence. We make a Farey symbol out of the generalized Farey
sequence by adding the pairing information in the obvious way.

On the other hand, if $F$ is a Farey symbol we can construct a
special polygon for $F$. For adjacent entries of the Farey sequence
$x_{i-1}$ and $x_i$, if $p_i$ is a free pairing or an even pairing
we let $P$ have as a side the hyperbolic arc joining $x_{i-1}$ and
$x_i$. Otherwise if it is odd we let $\gamma$ be the unique element
of $\PSLZ$ such that $\gamma(0) = x_{i-1}$ and $\gamma = x_i$ and
join $x_{i-1}$ and $x_i$ by the arcs $\gamma(H_{0,\rho})$ and
$\gamma(H_{\rho,\infty})$. Thus we get a hyperbolic polygon which is
made into a special polygon by adding pairing information in the
obvious way.

\begin{example}
$\Gamma(2)$ has a Farey symbol
 $\xymatrix@C=15pt{-\infty
 \ar@{-}@/_/[r]_1&
 \frac{0}{1} \ar@{-}@/_/[r]_{2} &
 \frac{1}{1} \ar@{-}@/_/[r]_{2} &
 \frac{2}{1} \ar@{-}@/_/[r]_1& \infty}$

\end{example}

\subsection{Generators}

If $P$ is a special polygon for a group $\Gamma$ then $\Gamma$ is
independently generated by the transformations mapping each side to
its paired side. If $F$ is a Farey symbol:
\begin{equation*}
 \xymatrix@C=15pt{-\infty
 \ar@{-}@/_/[r]_{p_0}&
 a_0/b_0 \ar@{-}@/_/[r]_{p_1} &
 a_1/b_1 \ar@{-}@/_/[r]_{p_2} &
 \dots \ar@{-}@/_/[r]_{p_{n-1}} &
 a_{n-1}/b_{n-1} \ar@{-}@/_/[r]_{p_n} &
 a_n/b_n \ar@{-}@/_/[r]_{p_{n+1}}& \infty}
\end{equation*}
then we can explicitly give formulas for the $\gamma$ corresponding
to a given side pairing.
\begin{theorem}
Suppose $(a_i/b_i, a_{i+1}/b_{i+1})$ are two adjacent vertices of
$F$. Then if the pairing between them $p_{i+1}$ is an even pairing,
let:
\begin{equation*}
G_{i+1} = \begin{pmatrix} a_{i+1}b_{i+1} + a_i b_i & -a_i^2 -
a_{i+1}^2 \\ b_i^2 + b_{i+1}^2 & -a_{i+1} b_{i+1} - a_i b_i
\end{pmatrix}
\end{equation*}
If $p_{i+1}$ is an odd pairing, let:
\begin{equation*}
G_{i+1} = \begin{pmatrix} a_{i+1}b_{i+1} + a_i b_{i+1} + a_i b_i &
-a_i^2 - a_i a_{i+1}-a_{i+1}^2 \\ b_i^2 + b_i b_{i+1} + b_{i+1}^2 &
-a_{i+1} b_{i+1} - a_{i+1} b_i - a_i b_i
\end{pmatrix}
\end{equation*}
And if $p_{i+1}$ is a free pairing that is paired with the side
between $a_k/b_k$ and $a_{k+1}/b_{k+1}$, let:
\begin{equation*}
G_{i+1} = \begin{pmatrix} a_{k+1}b_{i+1} + a_k b_i & -a_k a_i -
a_{k+1} a_{i+1} \\ b_k b_i + b_{k+1} b_{i+1} & -a_{i+1} b_{k+1} -
a_i b_k
\end{pmatrix}
\end{equation*}
Then $G_{i+1}$ is the side transformation corresponding to the
pairing $p_{i+1}$.
\end{theorem}
\begin{proof}
 \cite{kulkarni91} Theorem 6.1
\end{proof}

\subsection{Group Invariants} Several invariants of the group
$\Gamma$ can be read off from the Farey symbol $F$. Firstly, the
number of inequivalent order-2 (resp. order-3) elliptic points,
$e_2$ (resp. $e_3$), is the number of even (resp. odd) pairings in
$F$. Also, the number of free pairings in $F$ (half the number of
free edges) is equal to $r$, the rank of $\pi_1(\Gamma \backslash
\HH)$ (the fundamental group of the uncompactified modular curve).

To discuss the cusps of $\Gamma$, note that if $(x_i, x_{i+1})$ is
an edge with an even or odd pairing, then $x_i$ and $x_{i+1}$ are
equivalent cusps (since $G_{i+1} \in \Gamma$ maps $x_i$ to
$x_{i+1}$). Likewise, if $(x_i, x_{i+1})$ and $(x_j, x_{j+1})$ are
paired edges then $x_i$ and $x_{j+1}$ are equivalent cusps and $x_j$
and $x_{i+1}$ are equivalent cusps. This defines an equivalence
relation on the vertices of $P$. The equivalence classes are easy to
compute, because the defining equivalences occur in a cyclic patten.
So the number of cusps $t$ can be counted as the number of
equivalence classes.

For an edge $(\frac{a_i}{b_i}, \frac{a_{i+1}}{b_{i+1}})$ let $\gamma
= \begin{pmatrix}a_i & a_{i+1} \\ b_i & b_{i+1}\end{pmatrix}$. So
$\gamma^{-1}(x_i) = \infty$ and $\gamma^{-1}(x_{i+1}) = 0$. Then
define the \textbf{width} of a vertex $x_i$ to be the ``width'' of
$\gamma P$ at $\infty$. That is:
\begin{equation*}
 \text{width}(x_i) = \left\{
 \begin{array}{ll}
  |a_{i-1} b_{i+1} - a_{i+1} b_{i-1}| & \text{if $x_i$ is adjacent to no odd
  edge}\\
  |a_{i-1} b_{i+1} - a_{i+1} b_{i-1}| + 1/2 & \text{if $x_i$ is adjacent to 1 odd
  edge}\\
  |a_{i-1} b_{i+1} - a_{i+1} b_{i-1}| + 1& \text{if $x_i$ is adjacent to 2 odd
  edges}\\
 \end{array} \right.
\end{equation*}
The cusp width of a cusp $x$ of $\Gamma$ is then the sum of the
widths of the vertices of $P$ $\Gamma$-equivalent to $x$.

$\Gamma \backslash \HH$ is a genus $g$ orientable surface with $t$
points missing, one for each cusp. The rank of its fundamental group
is $r = 2g + t - 1$, so we can calculate the genus $g = \frac{r - t
+ 1}{2}$. Moreover, using the Hurwitz formula (\cite{shim1} Prop.
1.40) we get the index of $\Gamma$ in $\PSLZ$, $\mu = 3e_2 + 4e_3 +
12g + 6t - 12$. An even simpler formula for the index comes from
noting that $n + 2 = 2r + e_2 + e_3$ where $n+1$ is as in Definition
3. This, combined with the previous formula, implies  $\mu = 3n +
e_3$.

\section{Coset Permutation Representation of a Group}

Another method of representing groups that will be useful to us in
determining if a group is congruence is the coset permutation
representation developed by Millington \cite{Millington691}
\cite{Millington692}. Let $\Gamma$ be a subgroup of $\PSLZ$ with
$[\PSLZ: \Gamma] = \mu$ and $\PSLZ = \cup_{i=1}^\mu \alpha_i \Gamma$
a coset decomposition with $\alpha_1 = I$. Let $F$ be the standard
fundamental domain for $\PSLZ$. Then $\cup_{i=1}^\mu \alpha_i^{-1}
F$ is a fundamental domain for $\Gamma$. Let
\begin{equation*}
 E=\begin{pmatrix}0&1\\-1&0 \end{pmatrix}, \quad
 V=\begin{pmatrix}1&1\\-1&0 \end{pmatrix}, \quad
 L=\begin{pmatrix}1&1\\ 0&1 \end{pmatrix}, \quad
 R=\begin{pmatrix}1&0\\ 1&1 \end{pmatrix}
\end{equation*}

$E$ and $V$ generate $\PSLZ$, as do $L$ and $R$. The conversions
between them are:
\begin{eqnarray} \label{eq:lrdict}
 E = LR^{-1}L, \quad V = R^{-1}L\\
 L = EV^{-1}, \quad R = EV^{-2}
\end{eqnarray}

We have $E^2 = V^3 = 1$. In fact it is well-known that $\PSLZ$ is
isomorphic to the group (cf. \cite{Rankin-book-mf}):
\begin{equation} \label{eq:freegroup}
 \PSLZ \cong \langle e, v: e^2 = v^3 = 1 \rangle
\end{equation}

For each $\gamma$ in $\PSLZ$, left multiplication acts on the left
cosets of $\Gamma$ in $\PSLZ$ by permutation, i.e. there is a
homomorphism $\phi: \PSLZ \rightarrow S_n$, such that if
$\phi(\gamma) = \sigma_\gamma$ then $\gamma\alpha_i \Gamma =
\alpha_{\sigma_\gamma(i)} \Gamma$. In this way every finite-index
subgroup of $\PSLZ$ is associated with a pair of permutations
$e=\varphi(E)$ and $v=\varphi(V)$ with $e^2 = v^3 = 1$ which
generate a transitive permutation group (transitivity comes from $E$
and $V$ generating $\PSLZ$). We call $(e,v)$ a \textbf{coset
permutation representation} of $\Gamma$ and $(l,r)$ an
\textbf{LR-representation} of $\Gamma$, where $l = \varphi(L)$ and
$r = \varphi(R)$. Each form can be obtained from the other form by
the equations (\ref{eq:lrdict}). Note that $\gamma \in \PSLZ$ is in
$\Gamma$ if and only if $\gamma \Gamma = \Gamma$, i.e.
$\sigma_\gamma(1) = 1$.

On the other hand, suppose $e$ and $v$ are a pair of permutations on
$\mu$ letters with $e^2 = v^3 = 1$ that generate a transitive
permutation group $S$ (such a permutation we call {\textbf{valid}).
Define a homomorphism $\varphi: \PSLZ \rightarrow S$ such that
$\varphi(E) = e$ and $\varphi(V) = v$ (This is well-defined because
of (\ref{eq:freegroup})). Let $\Gamma = \{\gamma \in \PSLZ :
\varphi(\gamma)(1) = 1\}$. Then $\Gamma$ is an index-$\mu$ subgroup
of $\PSLZ$. Thus we have a correlation between valid pairs of
permutations and finite-index subgroups of $\PSLZ$. To test if $A
\in \PSLZ$ is in $\Gamma$ we write $A$ as a word in $L$ and $R$
(Using, essentially, the Euclidean Algorithm) and replace $L$ and
$R$ with the permutations $l$ and $r$. If the resulting permutation
fixes $1$ then $A$ is in $\Gamma$.

If one of the cosets is fixed by $e$, say $e(i) = i$, it corresponds
to an elliptic element in $\Gamma$, for $E \alpha_i \Gamma =
\alpha_i \Gamma$ means $\alpha_i^{-1} E \alpha_i \Gamma = \Gamma$,
meaning $\alpha_i^{-1} E \alpha_i$ (which is order 2) is in
$\Gamma$. So $e_2$, the number of inequivalent elliptic elements of
order 2 in $\Gamma$, is equal to the number of elements fixed by
$e$. Similarly, $e_3$ is the number of elements fixed by $v$.

The cusp width of $\Gamma$ at $\infty$ is the smallest positive
integer $n$ such that $L^n \in \Gamma$. Thus the cusp width at
infinity is the order of the cycle in $\varphi(L)$ which contains
``1''. Likewise, suppose $i$ is in a cycle of length $k$ in
$\varphi(L)$, i.e. $L^k \alpha_i \Gamma = \alpha_i \Gamma$, but $L^n
\alpha_i \Gamma \neq \alpha_i \Gamma$ for $0 < n < k$. Then
$\alpha_i^{-1} L^k \alpha_i \in \Gamma$, but $\alpha_i^{-1} L^n
\alpha_i \notin \Gamma$ for $0 < n < k$. If $q = \alpha_i^{-1}
\infty$ then $\alpha_i^{-1} L^k \alpha_i q = q$ but $\alpha_i^{-1}
L^n \alpha_i q \neq q$ for $0 < n < k$. Thus $\alpha_i^{-1} L^k
\alpha_i$ is a generator for the stabilizer of the cusp $q$, and
this cusp has width $k$.

\section{Algorithms}

\subsection{Calculating a Farey Symbol}
Recall that $T$ is the standard fundamental domain for $\PSLZ$, and
let $T^*$ be the hyperbolic triangle with vertices $\rho$, $i$ and
$\infty$ (So $T = T^* \cup (-\overline{T^*})$). $\mathscr{T} =
\{\gamma T : \gamma \in \PSLZ \}$ is a tessellation of the upper
half plane and any finite index subgroup $\Gamma$ has a fundamental
domain which is a simply connected union of $\mathscr{T}$-tiles. Let
$\mathscr{T}^* = \{\gamma T^* : \gamma \in \PSLZ \} \cup \{\gamma
(-\overline{T^*}) : \gamma \in \PSLZ \}$. $\mathscr{T}^*$ is also a
tessellation of the upper half plane, and we will construct a
fundamental domain for $\Gamma$ out of $\mathscr{T}^*$-tiles. The
starting point for our construction will be the six tiles around an
odd vertex. The following lemma shows this is a reasonable starting
point:

\begin{lemma}
Let $\Gamma$ be a subgroup of $\PSLZ$ with index $\geq 3$. Then the
stabilizer of $\rho = \frac12 + \frac{\sqrt{3}}{2} i$ or $\rho-1 =
-\frac12 + \frac{\sqrt{3}}{2} i$ is trivial (i.e. one of these
points is not elliptic in $\Gamma$).
\end{lemma}
\begin{proof}
If the two stabilizers are not trivial then they must be
$\Gamma_\rho = \left\{I,A, A^2 \right\}$ and $\Gamma_{\rho-1} =
\left\{I,B,B^2 \right\}$ where $A =
\begin{pmatrix}0&1\\-1&-1\end{pmatrix}$ and $B =
\begin{pmatrix}-1&1\\-1&0\end{pmatrix}$. But $A$ and $B$ generate
an index-2 subgroup of $\PSLZ$. So $\Gamma$ is either the (unique)
index-2 subgroup of $\PSLZ$ or $\PSLZ$ itself. And if the index of
$\Gamma$ in $\PSLZ$ is bigger than 2, at least one of $A$ and $B$
cannot be in $\Gamma$.
\end{proof}

\begin{figure}[t]\label{fig:5}
 \includegraphics*[angle=270,scale=.3, viewport=50 120 400 530]{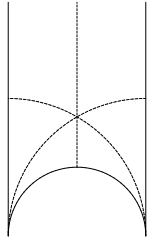}
 \caption{A hyperbolic triangle}
\end{figure}

So if $\Gamma$ is not $\PSLZ$ or $\Gamma_2$, the unique index 2
subgroup of $\PSLZ$ (cf. \cite{Rankin-book-mf}), then the hyperbolic
triangle with vertices either $0$, $1$ and $\infty$, or $-1$, $0$
and $\infty$ (cf. Figure 4) is contained in a fundamental domain of
$\Gamma$. The triangle is made of 6 $\mathscr{T}^*$-tiles. We will
make a polygon $P$ starting with this triangle, then attach
$\mathscr{T}^*$-tiles to $P$ and assign partial pairing information
to sides until we get a fundamental domain for $\Gamma$ (at which
point all the pairing information will be filled in). In the
algorithm we will say a $\mathscr{T}^*$-tile $T$ is
\textbf{adjoinable} to $P$ if $T$ is adjacent to a tile of $P$ and
if $P \cup T$ is contained in some fundamental domain of $\Gamma$.
Note that if $T$ is adjacent to $P$ with adjacency edge $e$ and if
$e$ cannot be paired with any other edge of $P$ then $T$ is
adjoinable.

Algorithm:

\begin{enumerate}
\item If $\Gamma = \PSLZ$ let $P$ be the special polygon with Farey
symbol
$$\xymatrix@C=15pt{-\infty \ar@{-}@/_/[r]_\circ&
 0 \ar@{-}@/_/[r]_\bullet& \infty}$$
or if $\Gamma = \Gamma_2$ let $P$ be the special polygon with Farey
symbol
$$\xymatrix@C=15pt{-\infty \ar@{-}@/_/[r]_\bullet&
 0 \ar@{-}@/_/[r]_\bullet& \infty}$$
In either case return $P$ and terminate.
\item If $\begin{pmatrix}-1&1\\-1&0\end{pmatrix}$ is not in $\Gamma$
then let $P$ be the hyperbolic polygon with vertices $0$, $1$, and
$\infty$. Otherwise let $P$ be the hyperbolic polygon with edges
$-1$, $0$ and $\infty$.
\item If any of the three sides of $P$ map to each other by a
$\gamma \in \Gamma$, assign that pairing to the side. (Note that
initially all sides are even sides).
 \item $P$ is now a polygon where every side is either:
 \begin{enumerate}
  \item even and already paired.
  \item odd and already paired.
  \item even and unpaired.
 \end{enumerate}
\item Pick an unpaired even side $e$. Figure 5 shows the
typical case (The other cases are the same as this case with
everything translated by some $\gamma \in \PSLZ$). Since $e$ is
unpaired, $T_1$ and $T_2$ must be adjoinable. If $o_1$ and $o_2$ are
the new odd edges of $P$ after adding $T_1$ and $T_2$ to $P$ then
either $\gamma o_1 = o_2$ for some $\gamma \in \Gamma$, or there is
no such $\gamma$. If there is $\gamma$ pair the two edges and go to
Step 3.

\begin{figure}[t]
 \includegraphics*[angle=270,scale=.5, viewport=50 120 450 530]{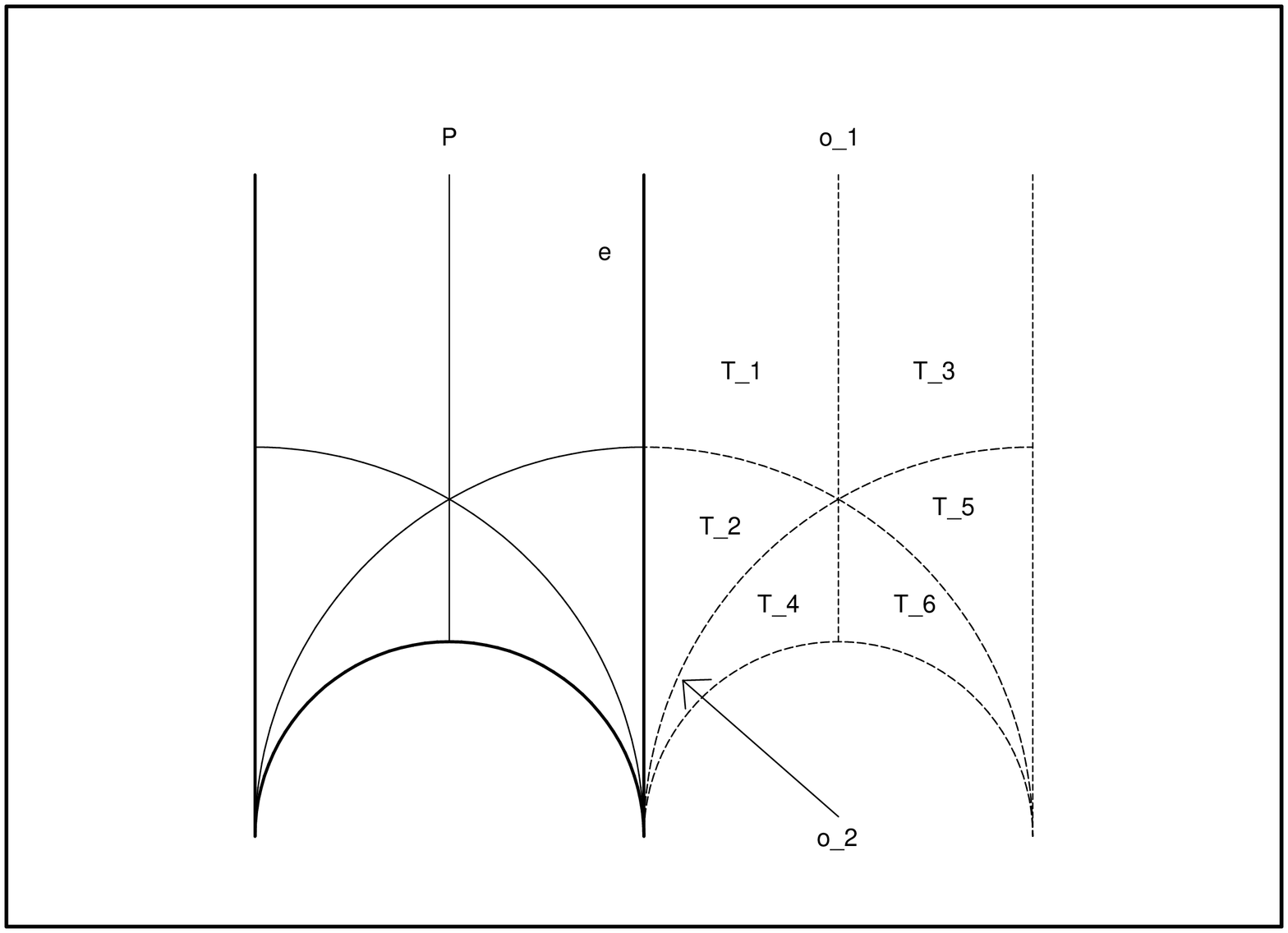}
 \caption{}
\end{figure}

\item If $o_1$ doesn't pair with $o_2$ then it doesn't pair with
any other side because the only other unpaired sides are odd. So
tiles $T_3$ and likewise $T_4$ are adjoinable. Each of these tiles
has a free edge and the free edges cannot pair with each other
(because their common vertex would have an internal angle of
$\frac{4\pi}{3}$, so the pairing transformations would make things
overlap), so $T_5$ and $T_6$ are adjoinable.
\item We've now added 6 $\mathscr{T}^*$-tiles to $P$ (One even
triangle). If either of the new even edges pair with any of the old
unpaired even edges then assign that pairing.
\item If all the sides of $P$ are paired then we are done.
Otherwise go to Step 4.
\end{enumerate}

The output of the algorithm is a special polygon $P$ with $\Gamma_P
= \Gamma$. Note that the algorithm must terminate, because a
fundamental domain of $\Gamma$ has hyperbolic area $\frac{\pi}{3}
[\PSLZ:\Gamma]$ and a single $\mathscr{T}^*$-tile has area
$\frac{\pi}{6}$. So for $P$ to be contained in a fundamental domain
of $\Gamma$ it can have at most $2 \cdot [\PSLZ:\Gamma]$
$\mathscr{T}^*$-tiles.

To effectively implement the algorithm we use Farey symbols. We need
only a way to test for group membership. Note that if $p_i/q_i$ and
$p_{i+1}/q_{i+1}$ are two adjacent vertices of the fundamental
polygon then the hyperbolic triangle added to the edge
$H_{x_i,x_{i+1}}$ in Step 6 is the triangle with vertices $p_i/q_i$,
$p_{i+1}/q_{i+1}$, and $(p_i+p_{i+1})/(q_i+q_{i+1})$.

So given a finite-index subgroup $\Gamma$ of $\PSLZ$, if we have a
way to test for group membership we can calculate a Farey symbol by
the following algorithm:

Algorithm for calculating a Farey Symbol:
\begin{enumerate}
\item If $\begin{pmatrix}1&1\\0&1 \end{pmatrix}$ and
$\begin{pmatrix}0&-1\\1&0 \end{pmatrix}$ are in $\Gamma$ then
$\Gamma = \PSLZ$, so return

$$\xymatrix@C=15pt{-\infty \ar@{-}@/_/[r]_\circ&
 0 \ar@{-}@/_/[r]_\bullet& \infty}$$

and terminate. If $\begin{pmatrix}0&1\\-1&-1 \end{pmatrix}$ and
$\begin{pmatrix}-1&1\\-1&0 \end{pmatrix}$ are in $\Gamma$ then
$\Gamma = \Gamma_2$, so return

$$\xymatrix@C=15pt{-\infty \ar@{-}@/_/[r]_\bullet&
 0 \ar@{-}@/_/[r]_\bullet& \infty}$$

and terminate.

\item If $\begin{pmatrix}-1&1\\-1&0 \end{pmatrix} \not\in \Gamma$ then
let $F$ be the (partial) Farey symbol: $$\xymatrix@C=15pt{-\infty
 \ar@{-}@/_/[r]_~&
 \frac{0}{1} \ar@{-}@/_/[r]_~ &
 \frac{1}{1} \ar@{-}@/_/[r]_~& \infty}$$
 Otherwise let $F$ be:
 $$\xymatrix@C=15pt{-\infty
 \ar@{-}@/_/[r]_~&
 \frac{-1}{1} \ar@{-}@/_/[r]_~ &
 \frac{0}{1} \ar@{-}@/_/[r]_~& \infty}$$

\item For each $i$ with $0\leq i \leq n+1$, if the pairing between
$x_{i-1}$ and $x_{i}$ is not filled in then check if it can be
paired with itself (even or odd pairing), or if it can be paired
with another unpaired edge (i.e., check if the appropriate $G_i$ is
in $\Gamma$). Wherever something can be paired, assign that pairing.

\item If all edges are now paired, return $F$ and terminate.

\item If there is still an unpaired edge, say between $p_i/q_i$ and
$p_{i+1}/q_{i+1}$, make a new vertex $(p_i+p_{i+1})/(q_i+q_{i+1})$
with no pairing information on the edges adjacent to it. Go to Step
3.
\end{enumerate}

The output is a Farey symbol for $\Gamma$.

\subsection{Group Membership}
The following algorithm described in \cite{llt95} tests if $A \in
\PSLZ$ is an element of the group corresponding to a Farey symbol
$F$. We will need a lemma about even lines:

\begin{lemma}
Let $l$ be an even line (a semicircle on the upper half plane with
rational endpoints $a/b$ and $a'/b'$ such that $|ab'-a'b| = 1$). Let
$P$ be a special polygon in $\HH$. Then either $l \subset P$ or $l
\cap P = \emptyset$.
\end{lemma}
\begin{proof}
\cite{llt95} Proposition 2.1
\end{proof}

Let $\Gamma$ be a finite index subgroup of $\PSLZ$ and $A$ an
element of $\PSLZ$.
\begin{equation}
A = \begin{pmatrix}c_0&c_0'\\d_0&d_0' \end{pmatrix}
\end{equation}
$A$ maps the even line $H_{0,\infty}$ to $l=H_{c'_0/d'_0,c_0/d_0}$.
By the lemma, either $l \subset P$ or it is disjoint from $P$
(except possibly at endpoints). If it is disjoint there is an edge
which it is naturally ``closest'' to (In a sense discussed in
\cite{llt95}). The idea of the algorithm is to translate $P$ across
the ``closest'' edge until $P$ intersects $H_{c/d,c'/d'}$, at which
point $A$ will be in $\Gamma$ if and only if $l$ is the image of
$(0,\infty)$ or the an edge paired with $(0,\infty)$. In the actual
algorithm we work in the other direction, translating the even line
instead of the special polygon.

Algorithm: \cite{llt95}

Let $k = 0$ and $F$ be a Farey symbol for $\Gamma$ with $0$ as one
of its vertices. Without loss of generality, we can assume
$\frac{c'_k}{d'_k} < \frac{c_k}{d_k}$.
\begin{enumerate}
\item There are two possibilities: If $\frac{c'_k}{d'_k}$ and
$\frac{c_k}{d_k}$ are both vertices of $P$ then terminate. Otherwise
we must have $x_i \leq \frac{c'_k}{d'_k} < \frac{c_k}{d_k} \leq
x_{i+1}$ with at least one ``$\leq$'' a strict inequality.

\item Let $g_{i+1}$ be the generator corresponding to the pairing
$p_{i+1}$ (recall this is the transformation mapping
$l=H_{c'_k/d'_k,c_k/d_k}$ to its paired side). If $p_{i+1}$ is a
free or even pairing, let $\alpha_k = g_{i+1}$. If $p_{i+1}$ is an
odd pairing, let $m = \frac{a_i+a_{i+1}}{b_i+b_{i+1}}$ where $x_i =
\frac{a_i}{b_i}$, $x_{i+1} = \frac{b_{i+1}}{b_{i+1}}$. Then the
interval $(\frac{c'_k}{d'_k},\frac{c_k}{d_k})$ must be between
either $x_i$ and $m$ or between $m$ and $x_{i+1}$. If
$\frac{c_k}{d_k} \leq m$, let $\alpha_k = g_{i+1}$. Otherwise let
$\alpha_k = g_{i+1}^{-1}$.

\item Let $\frac{c_{k+1}}{d_{k+1}} = \alpha_k \cdot
\frac{c_k}{d_k}$, $\frac{c'_{k+1}}{d'_{k+1}} = \alpha_k \cdot
\frac{c'_k}{d'_k}$. Replace $k$ with $k+1$ and go to Step 1.
\end{enumerate}

The algorithm returns $\frac{c'_k}{d'_k}$ and $\frac{c_k}{d_k}$,
which are two vertices of $P$, and a list of $\alpha_i$'s.
\begin{theorem}
The algorithm terminates, and $A$ is in $\Gamma$ if and only if one
of the following is true:
\begin{enumerate}
\item $\begin{pmatrix}c_k&c'_k \\d_k&d'_k \end{pmatrix} = \pm
\begin{pmatrix}1&0\\0&1 \end{pmatrix}$
\item $(\frac{c'_k}{d'_k}, \frac{c_k}{d_k})$ is a free side paired
with $(0,\infty)$.
\item $\begin{pmatrix}c_k&c'_k \\d_k&d'_k \end{pmatrix} = \pm
\begin{pmatrix}0&-1\\1&0\end{pmatrix}$ and $0$ and $\infty$ are
adjacent vertices with an even pairing between them.
\end{enumerate}
\end{theorem}
\begin{proof}
\cite{llt95}
\end{proof}
In addition, if $A$ is in $\Gamma$, $A$ can be written as a word in
the generators of $\Gamma$ because $A = \alpha_0^{-1} \alpha_1^{-1}
\dots \alpha_k^{-1} \begin{pmatrix}c_k&c'_k\\d_k&d'_k\end{pmatrix}$,
and each term in that product is one of the generators for $F$.

\subsection{Coset Representatives}
Let $\Gamma$ be a group with special polygon $P$. Let $T$ be the
hyperbolic triangle with vertices $i$, $\rho$, and $\infty$. By the
construction of $P$, $T$ is contained in $P$. The set of $\gamma \in
\Gamma$ such that $\gamma T$ is in $P$ is a set of coset
representatives of $\Gamma$.

Let $a_i/b_i$ and $a_{i+1}/b_{i+1}$ be a vertex of the special
polygon, and let $T =
\begin{pmatrix}1&1\\0&1\end{pmatrix}$ and $\varphi =
\begin{pmatrix}a_i & a_{i+1} \\ b_i & b_{i+1} \end{pmatrix}$. Then
$\varphi^{-1}(a_i/b_i) = \infty$ and $\varphi^{-1}(a_{i+1}/b_{i+1})
= 0$. Let $w_i$ be $|a_{i-1} b_{i+1} - a_{i+1} b_{i-1}|$ if the
pairing between $a_i/b_i$ and $a_{i+1}/b_{i+1}$ is not an odd
pairing and $|a_{i-1} b_{i+1} - a_{i+1} b_{i-1}| + 1$ if it is. Then
$w_i$ is the number of $\mathscr{T}^*$-tiles of the form $\gamma T$
in $P$. Thus a list of left coset representatives for $\Gamma$ is
$\bigcup_{i=0}^n \{T^{-j} \phi_i^{-1} : 0 \leq j < w_i \}$.

\subsection{Congruence Testing}
Let $\Gamma$ be a finite index subgroup of $\PSLZ$. Lang, Lim and
Tan give a test purely in terms of Farey symbols to determine if
$\Gamma$ is a congruence group \cite{llt95}. Their test relies on
Wohlfahrt's Theorem \cite{wohlfahrt64} which says that if $\Gamma$
has level $N$ then $\Gamma$ is a congruence group if and only if
$\Gamma$ contains $\Gamma(N)$. In Lang, Lim and Tan's test, if
$\Gamma$ has level $N$ one computes a Farey symbol for $\Gamma(N)$,
giving a complete set of generators for $\Gamma(N)$. One then checks
if each of these generators is contained in $\Gamma$ using the above
algorithm. The difficulty with this algorithm is that the index of
$\Gamma(N)$ increases very quickly with $N$, so if $\Gamma$ has
large level, the calculation of a Farey symbol for $\Gamma(N)$ can
be very lengthy, even if $\Gamma$ has relatively small index.

Another test for congruence was developed by Tim Hsu using
Millington's coset permutation representations \cite{hsu96}. If we
have an LR-representation of $\Gamma$ there is a list of relations
that are satisfied if and only $\Gamma$ is congruence.

To calculate an LR-representation from a Farey symbol, use the above
algorithm to calculate a list of left coset representatives
$\alpha_i \in \PSLZ$ where $\PSLZ = \bigcup_{i=1}^\mu \alpha_i
\Gamma$. To calculate $l$, for instance, recall that $l$ is the
permutation such that $L\alpha_i \Gamma = \alpha_{l(i)} \Gamma$. So
$l$ sends $i$ to the unique $j$ such that $\alpha^{-1}_j L \alpha_i
\in \Gamma$. So we run through every $1\leq i \leq \mu$ and
calculate the permutation. $r$ can be calculated similarly.
(Actually, although we need $l$ and $r$ it is easier to calculate
$e$ and $v$, because we know beforehand that they are order 2 and 3
respectively. Then $l = ev^{-1}$ and $r = ev^{-2}$).

Knowing $l$ and $r$ we can directly apply Tim Hsu's congruence
algorithm \cite{hsu96}. Depending on the order of $l$, (i.e. the
level of $\Gamma$) there are different lists of relations of $l$ and
$r$ that are satisfied if and only if $\Gamma$ is congruence. For
example, if $N$ is the order of $l$ and $N$ is odd then $\Gamma$ is
a congruence group if and only if $r^2l^{-\frac{1}{2}}$ is the
identity permutation (where $\frac{1}{2}$ is the inverse of $2$
modulo $N$).

\section{Implementation}

Helena Verrill has a written a MAGMA package for working with Farey
symbols for congruence groups. Also, for congruence or noncongruence
groups, the algorithms described above have been implemented by the
first author as a collection of functions for SAGE. The package and
basic examples may be downloaded at:

\url{http://www.public.iastate.edu/~kurthc/research/index.html}

%\bibliography{kurth}{}
%\bibliographystyle{amsalpha}
\providecommand{\bysame}{\leavevmode\hbox
to3em{\hrulefill}\thinspace}
\providecommand{\MR}{\relax\ifhmode\unskip\space\fi MR }
\providecommand{\MRhref}[2]{%
  \href{http://www.ams.org/mathscinet-getitem?mr=#1}{#2}
} \providecommand{\href}[2]{#2}

%\bibliographystyle{amsalpha}
%\bibliography{longbibl}

\end{document}